\documentclass{amsart}
\usepackage{amssymb}
\usepackage{amsthm}
\usepackage{amsmath}
\usepackage[curve,matrix,arrow]{xy}
\theoremstyle{plain}\newtheorem{Theorem}{Theorem}[section]
\theoremstyle{plain}
\theoremstyle{plain}\newtheorem{Corollary}[Theorem]{Corollary}
\theoremstyle{plain}\newtheorem{Lemma}[Theorem]{Lemma}
\theoremstyle{plain}\newtheorem{Proposition}[Theorem]{Proposition}
\theoremstyle{plain}
\theoremstyle{definition}\newtheorem{Definition}[Theorem]{Definition}
\theoremstyle{definition}\newtheorem{Example}[Theorem]{Example}
\theoremstyle{definition}
\theoremstyle{definition}\newtheorem{Remark}[Theorem]{Remark}
\theoremstyle{definition}
\theoremstyle{definition}

  \def\OG{{\mathcal{O}G}}

\def\CH{{\mathcal{H}}}

\def\CO{{\mathcal{O}}}

\def\Q{{\mathbb Q}}
\def\Z{{\mathbb Z}}

\def\Br{\mathrm{Br}}             

\def\dim{\mathrm{dim}}           
\def\End{\mathrm{End}}           
\def\Endbar{\underline{\mathrm{End}}}
           
   \def\hatExt{\widehat{\mathrm{Ext}}}

\def\Hom{\mathrm{Hom}}           
\def\Hombar{\underline{\mathrm{Hom}}}

\def\Id{\mathrm{Id}}

\def\Irr{\mathrm{Irr}}           
           \def\tenO{\otimes_{\mathcal{O}}}
\def\mod{\mathrm{mod}}
\def\modbar{\underline{\mathrm{mod}}}
\def\op{\mathrm{op}}

\def\pr{\mathrm{pr}}
           
\def\rank{\mathrm{rank}}         
\def\rk{\mathrm{rk}}

\def\soc{\mathrm{soc}}

\def\Tr{\mathrm{Tr}}             
\def\tr{\mathrm{tr}}

\def\ualpha{\underline{\alpha}}

\newcommand{\inv}{^{-1}}

\title{On Tate duality and a projective scalar property for symmetric 
algebras} 
\author{Florian Eisele}
\address{City University London\\ Northampton Square\\ London\\ EC1V 0HB\\ UK}
\email{Florian.Eisele@city.ac.uk}

\author{Michael Geline}
\address{Northern Illinois University\\ DeKalb, IL 60115, USA}
\email{mgeline@niu.edu}

\author{Radha Kessar}
\address{City University London\\ Northampton Square\\ London\\ EC1V 0HB\\ UK}
\email{Radha.Kessar.1@city.ac.uk}

\author{Markus Linckelmann}
\address{City University London\\ Northampton Square\\ London\\ EC1V 0HB\\ UK}
\email{Markus.Linckelmann.1@city.ac.uk}


\begin{document}

\maketitle

\begin{abstract}
We identify a class of symmetric algebras over a complete discrete 
valuation ring $\CO$ of characteristic zero to which the 
characterisation of Kn\"orr lattices in terms of stable endomorphism 
rings in the case of finite group algebras, can be extended.  
This class includes finite group algebras, their  blocks and source algebras and  
Hopf orders.  We  also show that certain  arithmetic properties of finite group representations 
extend to  this class of algebras. Our 
results are based on an  explicit description of Tate duality for 
lattices over symmetric $\CO$-algebras whose extension to the quotient 
field  of $\CO$  is separable.
\end{abstract}

%%%%%%%%%%%%%%%%%%%%%%%%%%%%%%%%%%%%%%%%%%%%%%%%%%%%%%%%%%%%%%%%%%%%%%%%
\section{Introduction}

Let $p$ be  a prime. Let $\CO$ be a complete discrete 
valuation ring with maximal ideal $J(\CO)=$ $\pi\CO$ for some
$\pi\in$ $\CO$, residue field $k=$ $\CO/J(\CO)$ of characteristic 
$p$, and field of fractions $K$ of characteristic zero. 
An $\CO$-algebra $A$ is {\it symmetric} if $A$ is isomorphic to
its $\CO$-dual $A^*$ as an $A$-$A$-bimodule; this implies that $A$ is 
free  of finite rank over $\CO$. 
The image $s$ of $1_A$ under a bimodule isomorphism
$A\cong$ $A^*$ is called a {\it symmetrising form} for $A$; it
has the property that $s(ab)=$ $s(ba)$ for all $a, b\in A$ and
that the bimodule isomorphism $A\cong$ $A^*$ sends $a\in A$ to
the map $s_a\in$ $A^*$ defined by $s_a(b) = s(ab)$ for all
$a, b\in A$. Since the automorphism group of $A$ as an
$A$-$A$-bimodule is canonically isomorphic to $Z(A)^\times$,
any other symmetrising form of $A$ is of the form $s_z$ for
some $z\in$ $Z(A)^\times$. If $X$ is an $\CO$-basis of $A$, then any
symmetrising form $s$ of $A$ determines a dual basis $X^\vee=$
$\{x^\vee\ |\ x\in X\}$ satisfying $s(xx^\vee)=1$ for $x\in X$
and $s(xy^\vee)=0$ for $x, y\in$ $X$, $x\neq y$. We denote
by $\Tr^A_1 : A\to$ $Z(A)$ the $Z(A)$-linear map defined by
$\Tr^A_1(a)=$ $\sum_{x\in X} xax^\vee$ for all $a\in$ $A$. This
map depends on the choice of $s$  but not on the choice of the
basis $X$. We set $z_A=$ $\Tr^A_1(1_A)$ and call $z_A$ the
{\it relative projective element of $A$ in $Z(A)$ with respect to
$s$}. This is also called the {\it central Casimir element} in
\cite{Br09}.
If $z\in$ $Z(A)^\times$ and $s'=$ $s_z$, then the dual basis
of $X$ with respect to $s'$ is equal to $X^\vee z^{-1}$, where 
$X^\vee$ is the dual basis of $X$ with respect to $s$, and hence
the relatively projective element in $Z(A)$ with respect
to $s'$ is equal to $z'_A=$ $z_Az^{-1}$. If we do not specify a
symmetric form of a symmetric algebra $A$, then the relative 
projective elements form a $Z(A)^\times$-orbit in $Z(A)$. See Brou\'e \cite{Br09} for more details. 

The purpose of this paper is to examine situations in which some relative projective element is a scalar multiple of the identity.

\begin{Definition}
A symmetric $\CO$-algebra $A$ is said to have the {\it projective
scalar property} if there exists a symmetrising form $s$ of $A$
such that the corresponding relative projective element $z_A$ is of the 
form $z_A =$ $\lambda 1_A$ for some $\lambda \in$ $\CO$.
\end{Definition}

Throughout the paper we will be working with a symmetric 
$\CO$-algebra  $A$ such that  the $K$-algebra $K \tenO A$ is  
separable. Since $K$ has characteristic zero, $K\tenO A$ is separable  
if and only if it is semisimple. This in turn  is equivalent  
to the condition that the relative projective element with respect to 
some, and hence any, symmetrising form on $A$  is invertible in 
$Z(K \tenO A) $  (see \cite[Proposition 3.6]{Br09}). In  particular,  
in case $A$ has the  projective scalar property, the separability of 
$K \tenO A$ is equivalent to the property that  the relative 
projective elements of $A$ are non-zero.

Matrix algebras, finite group algebras, blocks  and source algebras   
of finite group algebras, as well as  Hopf algebras  whose extension 
to  $K$ is semisimple  have the projective scalar property (see 
Examples \ref{basicexample},  \ref{sourcealgebraexample}, and 
\ref{hopfexample}), but Iwahori-Hecke algebras and rings of 
generalized characters do not typically have this property (see 
Examples \ref{rank2example}, \ref{heckeexample},  and 
\ref{characterringexample}).  The projective scalar property is 
invariant under taking direct factors and tensor products but not  
under direct products, and is not invariant under Morita equivalences 
(see Example \ref{basicexample}).

Our motivation for studying  algebras with the projective scalar
property  comes  from   a  characterisation of   Kn\"orr  lattices  
for a finite group algebra  in terms of  the relatively $\CO$-stable  
module category of the algebra.  Recall that  an $A$-lattice is a  
left unital $A$-module which is free of finite rank as an $\CO$-module. 
An indecomposable $A$-lattice $U$ is called a {\it Kn\"orr lattice}
if the linear trace form $\tr_U$ on $\End_\CO(U)$ satisfies 
$\tr_U(\alpha)\CO\subseteq$ $\rk_\CO(U)\CO$ for every $\alpha \in \End_A(U)$, with equality precisely when $\alpha$ is an automorphism.

 Now for two finitely
generated $A$-modules $U$ and $V,$ we denote by $\Hombar_A(U,V)$ the homomorphism 
space in the $\CO$-stable category $\modbar(A)$ of finitely generated 
$A$-modules; that is, $\Hombar_A(U,V)$ is the quotient of $\Hom_A(U,V)$
by the subspace $\Hom_A^\pr(U,V)$ of $A$-homomorphisms $U\to$ $V$
which factor through a relatively $\CO$-projective $A$-module.
We write $\End_A^\pr(U)=$ $\Hom_A^\pr(U,U)$ and $\Endbar_A(U)=$
$\Hombar_A(U,U)$. 
 
For an $A$-lattice $U,$ let $a(U)$ denote the smallest non-negative 
integer such that $\pi^{a(U)}$ annihilates $\Endbar_A(U).$ In 
\cite{CarJo}, the element $\pi^{a(U)}$ is referred to as the {\it 
exponent} of $U$. If $U$ is indecomposable non-projective, $U$ is said 
to have the {\it stable exponent property} if the socle of 
$\Endbar_A(U)$ as a (left or right) module over itself is equal to 
$\pi^{a(U)-1}\Endbar_A(U).$  

%The stable exponent property is obviously 
%invariant under stable and in particular Morita equivalences.

%We shall see in Example \ref{knorrmoritaexample} that the property of being a Kn\"{o}rr lattice is not invariant under Morita equivalences. However, it is clear that the stable exponent property is invariant under such equivalences. 

Carlson and Jones \cite{CarJo}, and independently Thevenaz 
\cite{Thev88Dual} and Kn\"orr \cite{KnoUnpub} proved 
that for $G$ a finite group, an absolutely indecomposable 
non-projective $\OG$-lattice is a Kn\"orr 
lattice if and only if it has the stable exponent property. The projective scalar property guarantees such an equivalence:

%Since the property of being a
%Kn\"orr lattice is not invariant under Morita equivalences, the 
%Carlson-Jones-Kn\"orr-Thevenaz  result singles out algebras within 
%their Morita equivalence class for which the equivalence of these two 
%properties holds.   The projective scalar property  guarantees such 
%an equivalence.  

\begin{Theorem} \label{KnoerrExponent}
Let $A$ be a symmetric $\CO$-algebra such that $K\tenO A$ is separable.
Suppose that $A$ has the projective scalar property. 
Then an indecomposable non-projective $A$-lattice $U$ is a 
Kn\"orr lattice if and only if $U$ is absolutely indecomposable and
has the stable exponent property.
\end{Theorem}

The converse to this theorem is false. In Example \ref{non-converseexample}, we shall see a symmetric algebra without the projective scalar property for which the Kn\"{o}rr lattices coincide with those having the stable exponent property. Thus, the equivalence between the Kn\"{o}rr and stable exponent properties does not provide a characterization of the projective scalar property. Also, in Example \ref{knorrmoritaexample}, we shall see both Kn\"{o}rr lattices which do not have the stable exponent property, as well as lattices with the stable exponent property which are not Kn\"{o}rr. 

Example \ref{knorrmoritaexample} will, in addition, show that the property of being a Kn\"{o}rr lattice is not invariant under Morita equivalences. However, it is easy to see that the stable exponent property is invariant under such equivalences. Thus, two subclasses can be identified within a given Morita equivalence class of symmetric algebras: namely, those for which the above two types of lattices coincide, and those with the projective scalar property.  

%The converse to Theorem \ref{KnoerrExponent} does not hold in general. 
%Specifically, there are symmetric algebras without the projective 
%scalar property for which the Kn\"{o}rr lattices coincide with those 
%having the stable exponent property. See Example 
%\ref{non-converseexample}.  Also, in Example \ref{knorrmoritaexample}, 
%we shall see both Kn\"{o}rr lattices that do not have the stable 
%exponent property, as well as lattices with the stable exponent 
%property which are not Kn\"{o}rr. 

The basic ingredient for the proof of Theorem \ref{KnoerrExponent} 
is a description of Tate duality 
for lattices over symmetric $\CO$-algebras with separable coefficient 
extensions which makes the role of the relative projective element explicit.  Note that $\Hombar_A(U,V)$ is a torsion
$\CO$-module for any $A$-lattices $U$ and $V$ when $K\tenO A$ is separable. This follows from the 
Gasch\"utz-Ikeda Lemma (cf. \cite[Lemma 7.1.11]{GP}), which is a 
special case of Higman's criterion for modules over symmetric algebras 
in Brou\'e \cite{Br09}.

\begin{Theorem} \label{TatedualityO}
Let $A$ be a symmetric $\CO$-algebra with symmetrising form $s$
such that $K\tenO A$ is separable. Set $z=z_A$.
Let $U$, $V$ be $A$-lattices. The map sending
$(\alpha,\beta)\in$ $\Hom_A(U,V)\times \Hom_A(V,U)$ to
$\tr_{K\tenO U}(z^{-1} \beta\circ\alpha)\in$ $K$ induces a
non degenerate pairing
$$\Hombar_A(U,V)\times\Hombar_A(V,U)\to K/\CO\ .$$
\end{Theorem}

Here $\tr_{K\tenO U}(z^{-1} \beta\circ\alpha)$ is the trace of the
$K$-linear endomorphism of $K\tenO U$ obtained from extending
the endomorphism $\beta\circ\alpha$ of $U$ linearly to $K\tenO U$,
composed with the endomorphism given by multiplication on $K\tenO U$
with the inverse $z^{-1}$ of $z$ in $Z(K\tenO A)$.
If $A$ has the projective scalar 
property, then the Tate duality pairing admits the following 
description (which is in this form well-known for finite group
algebras; see \cite[Theorem (7.4)]{Brown}). 

\begin{Corollary} \label{TateDuality}
Let $A$ be a symmetric $\CO$-algebra such that $K\tenO A$ is 
separable. Suppose that $z_A=\pi^n 1_A$ for some choice of a
symmetrising form of $A$ and some positive integer $n$.
Let $U$ and $V$ be $A$-lattices. 
The map sending $(\alpha,\beta)\in$ $\Hom_A(U,V)\times \Hom_A(V,U)$ to 
$\tr_U(\beta\circ\alpha)$ induces a non degenerate pairing
$$\Hombar_A(U,V)\times\Hombar_A(V,U)\to \CO/\pi^n\CO\ .$$
\end{Corollary}

\begin{Remark}
Theorem \ref{TatedualityO}, applied to $U=V$, shows that if
$U$ is an indecomposable nonprojective lattice for a symmetric
$\CO$-algebra $A$ such that $K\tenO A$ is separable, then
the socle of $\Endbar_A(U)$ as a module over itself is simple, since
it is dual to $\Endbar_A(U)/J(\Endbar_A(U))\cong k$. This fact is
well-known - see Roggenkamp \cite{Rog77} - and this is the key step 
in the existence proof of almost split sequences of $A$-modules. 
Applying Theorem \ref{TatedualityO} to Heller translates of $V$
yields non degenerate pairings
$$\hatExt_A^n(U,V)\times\hatExt_A^{-n}(V,U) \to K/\CO$$
for any integer $n$. Applied to $U=V=A$ as a module over
$A\tenO A^\op$ this yields non degenerate pairings in
Tate-Hochschild cohomology
$$\widehat{HH}^n(A) \times \widehat{HH}^{-n}(A) \to K/\CO\ .$$
\end{Remark} 

\iffalse
The fact that being a Kn\"orr lattice is not invariant under
Morita equivalences raises the question as to how one would
need to modify the definition in such a way that it does become
invariant under Morita equivalences while still coinciding with
the notion of Kn\"orr lattices whenever the algebra has the
projective scalar property.   \fi

Theorem \ref{KnoerrExponent} is a special case of  the following 
consequence of Theorem~\ref{TatedualityO} which  gives a 
characterisation of  absolutely  indecomposable  modules  with the  
stable exponent property for symmetric  $\CO$-algebras.
 Denote by $\nu$ a $\pi$-adic valuation on $K$.

\begin{Theorem} \label{GenKnorr}
Let $A$ be a symmetric $\CO$-algebra with symmetrising form $s$
such that $K\tenO A$ is separable. Denote by $z$ the associated
relatively projective element of $A$ in $Z(A)$.
Let $U$ be an indecomposable non projective $A$-lattice.
The following are equivalent.

\begin{enumerate}
\item[\rm (i)]
For any $\alpha\in$ $\End_A(U)$ we have 
$\nu(\tr_{K\tenO U}(z^{-1}\alpha))\geq$ 
$\nu(\tr_{K\tenO U}(z^{-1}\Id_U))$, with equality if and only if
$\alpha$ is an automorphism of $U$.

\item[\rm (ii)]
The $A$-lattice $U$ is absolutely indecomposable and has the stable 
exponent property.
\end{enumerate}
\end{Theorem}

 Symmetric $\CO$-algebras with split semisimple coefficient extensions 
to $K$ having the projective scalar property can be characterised as 
follows.

\begin{Theorem} \label{scalarA}
Let $A$ be a symmetric $\CO$-algebra such that $K\tenO A$ is
split semisimple. Denote by $\rho : A\to$ $\CO$ the regular
character of $A$. The following are equivalent.

\begin{enumerate}
\item[\rm (i)]
The algebra $A$ has the projective scalar property.

\item[\rm (ii)]
There exists a non-negative integer $n$ such that
$\pi^{-n}\rho$ is a symmetrising form of $A$.

\item[\rm (iii)]
There exists a non-negative integer $n$ such that for any
$A$-lattice $U$ we have 
$$\tr_U(\End_A(U))= \pi^{n-a(U)}\CO\ .$$
\end{enumerate}

\noindent
Moreover, if these three equivalent statements hold, then the integers
$n$ in (ii) and (iii) coincide, and $\pi^n 1_A$ is a relative
projective element with respect to some symmetrising form of $A$.
\end{Theorem}

We also have a characterisation, in terms of the decomposition matrix,
  of symmetric $\CO$-algebras $A$ such that some algebra in the  
Morita  or derived  equivalence class of $A$ has the scalar   
projective property. Recall that if  $B$ is a  split finite 
dimensional algebra over a field $F$  then  the set of characters of 
simple $A$-modules is a linearly independent subset of  the $F$-vector 
space  of functions from $B$ to $F$ (see for instance 
\cite[Chapter 3, Theorem 3.13]{NT}), and hence may be  identified with 
a set of representatives of the isomorphism classes  of   simple 
$B$-modules.

 \begin{Theorem} \label{MoritascalarA}    
Let $A$ be a symmetric $\CO$-algebra such that $K\tenO A$ is split 
semisimple  and $ k\tenO A $ is split. Denote by $\Irr_K(A)$   the  
set of  characters of simple   $K \tenO  A$  modules  and by 
$\Irr_k(A)$  the set of characters  of  simple $ k\tenO A $-modules. 
For $\chi \in \Irr_K(A)$ and $\varphi\in \Irr_k(A)$ denote by 
$ d_{\chi, \varphi} $ the   multiplicity  of $S$ as a composition  
factor of  $k\tenO V $,  where  $V$ is an $A$-lattice  such that 
$ K \tenO   V$  has character $\chi $, and $S$ is  a simple 
$k\tenO A$-module  with character  $\varphi$. The following are 
equivalent.

\begin{enumerate}
\item[\rm (i)]
There exists an algebra Morita equivalent to $A$ with the projective 
scalar property.

\item[\rm (ii)]
There exists an algebra derived equivalent to $A$ with the  projective 
scalar property.

\item[\rm (iii)]
There exists a non-negative integer $n$ and positive integers  
$ m_{\varphi}  $, $\varphi \in \Irr_k(A)  $     such that    setting  
$ a_\chi:= \sum _{\varphi \in \Irr_k(A)} m_\varphi d_{\chi, \varphi} $,
$ \chi \in  \Irr_K(A) $, the form 
$\pi^{-n} \sum_{ \chi \in \Irr_K(A)} a_{\chi}  \chi $  is  a 
symmetrising form for $A$.

 \item[\rm (iv)] 
There exists a non-negative integer $n$  and integers $m_{\varphi}$, 
$\varphi \in \Irr_k(A) $  such that setting  
$a_{\chi}:= \sum _{\varphi \in \Irr_k(A)} m_\varphi d_{\chi, \varphi}$,
$ \chi \in \Irr_K(A) $, the form 
$\pi^{-n} \sum_{ \chi  \in \Irr_K(A)} a_{\chi} \chi$ is a symmetrising 
form for $A$.
\end{enumerate}

\end{Theorem}

 We point out that certain arithmetic features of finite group 
representations carry over to algebras with the projective scalar
property. Recall that the degree of an ordinary irreducible character 
of a finite group $G$ divides the order of $G$ and that if $U$ is a  
projective $\CO G$-lattice, then the $p$-part of $|G|$ divides the 
$p$-part of the $\CO$-rank of $U$.

\begin{Proposition} \label{scalardegree}
Let $A$ be a symmetric $\CO$-algebra such that $K\tenO A$ is split 
semisimple. Assume that $A$ has the projective scalar property  and 
let $\pi^n 1_A $  be a relative projective element    with respect to 
some symmetrizing form on $A$.

\begin{enumerate}
\item[\rm(i)]  If $U$ is a Kn\"{o}rr $A$-lattice, then the $p$-part of 
the $\CO$-rank of $U$ divides  $ \pi^n $ in $\CO.$

\item [\rm(ii)]   If $U$ is a projective  $A$-lattice, then  the 
$p$-part of the $\CO$-rank  of $U$ is  divisible in $\CO$ by $\pi^n$.
\end{enumerate}
\end{Proposition}

\begin{Remark} 
Note that if $A=\CO G$, then $|G| \cdot1_{\CO G} $ is the relative 
projective element with respect to the standard symmetrising form  
(see    \cite[Examples and Remarks after Proposition 3.3] {Br09}).  
Moreover,  an  absolutely irreducible  $\CO G$-lattice  is a 
Kn\"{o}rr $\CO G$-lattice.   Hence,  letting $p$  vary across   all 
primes in (i), one sees that  the above does generalise the 
corresponding results  for group algebras.   A  related  global  
divisbility  criterion for irreducible lattices of symmetric algebras
has been given by Jacoby and  Lorenz \cite[Corollary 6]{JaLo} in the context of Kaplansky's sixth conjecture.
\end{Remark}

For $U$ an $A$-lattice, define the height of $U$ to be the number 
$h(U)$ such that 
\[\rank(U)_p = p^{m+h(U)},  \]
where  $m$ is defined  by
\[
p^m = \text{min}_V\{ \rank(V)_p\}
\]    
as  $V$ ranges over all  irreducible $A$-lattices.  Note that $h(U)$ 
is a non-negative integer.
 
 It is well known that a Morita equivalence between  blocks of finite 
group algebras or between a block algebra and  the  corresponding 
source algebra preserves the height of corresponding irreducible   
characters (see \cite{BroueIsom} and \cite{BroueEq}).  The following 
theorem generalises this to  algebras with the projective scalar 
property and to   Kn\"{o}rr lattices.

\begin{Theorem}  \label{scalarheight}  
Suppose that $K \tenO A $ is split semisimple. Let  $A'$ be  an 
$\CO$-algebra  Morita equivalent to $A$, and suppose that both $A$ and 
$A'$  have  the projective scalar property. Let $U$ be a Kn\"{o}rr   
$A$-lattice and let  $U'$  be an $A'$-lattice corresponding to  $U$ 
through a Morita equivalence between $A$ and $A'$.  Then  $U'$  is a 
Kn\"{o}rr  $A'$-lattice  and  $h(U)= h(U')$.  
\end{Theorem}

Finally we point out that although the stable exponent property does 
not apply to projective lattices, we can, following Kn\"{o}rr 
\cite[Lemma 1.9]{Kno89}, characterise projective Kn\"{o}rr lattices 
in the presence of the projective scalar property.

\begin{Proposition}\label{knorrproj}
Let $A$ be as in the previous proposition. Assume that $U$ is an 
$A$-lattice which is both projective and Kn\"{o}rr. Then $U/\pi U$ is 
a simple $A/\pi A$-module. In particular, $K\tenO U$ is an irreducible 
$K\tenO A$-module.
\end{Proposition}
 
 Section \ref{Tate}  contains the proof of Theorem \ref{TatedualityO}, 
Theorem \ref{GenKnorr} and Theorem \ref{KnoerrExponent}. We prove 
Theorems \ref{scalarA} and \ref{MoritascalarA}  in Section 
\ref{sec_el_char}. This section also contains a  characterisation of 
the projective scalar property in terms of rational centres.
Section \ref{heights} discusses arithmetic properties of Kn\"{o}rr 
lattices  in the presence of the scalar projective property, including 
the proof of Proposition \ref{scalardegree}  and Theorem 
\ref{scalarheight}. Section \ref{examples} contains various examples.

%%%%%%%%%%%%%%%%%%%%%%%%%%%%%%%%%%%%%%%%%%%%%%%%%%%%%%%%%%%%%%%%%%%%%%%
\section{Tate Duality for Symmetric algebras }\label{Tate}

The proof of Theorem \ref{TatedualityO} is an adaptation of 
ideas in Th\'evenaz \cite[Section 1]{Thev88Dual}. 
We keep the notation in Theorem \ref{TatedualityO}. For simplicity,
we write in this section $KA=$ $K\tenO A$, $KU=$ $K\tenO U$, 
and $KV=$ $K\tenO V$.
We write $K\Hom_\CO(U,V)=$ $K\tenO \Hom_\CO(U,V)$ and identify this
space with $\Hom_{K}(KU,KV)$ whenever convenient. Similarly, 
we write $K\Hom_A(U,V)=$ $K\tenO \Hom_A(U,V)$ and identify this
space with $\Hom_{KA}(KU,KV)$.
Let $X$, $X^{\vee}$ be a pair of $\CO$-bases of $A$ dual to each other 
with respect to the symmetrising form $s$; in particular,
the relative projective element with respect to $s$ is
$$z_A= \sum_{x\in X}xx^\vee = \sum_{x\in X}x^\vee x\ ,$$ 
where $x^\vee$ denotes the unique element in $X^\vee$ satisfying 
$s(xx^\vee)=1$, for $x\in$ $X$. We denote by  
$$\Tr_1^A : K\Hom_{\CO}(U, V)\to K\Hom_{A}(U,V)$$ 
the $K$-linear map sending $\alpha\in$ $\Hom_{\CO}(U,V)$ 
to $\sum_{x\in X} x\alpha x^\vee$. Here $x\alpha x^\vee\in$
$\Hom_\CO(U,V)$ is defined by $(x\alpha x^\vee)(u)=$
$x\alpha(x^{\vee} u)$ for $u\in$ $U$ and $x\in$ $X$.
Clearly, $\Tr_1^A$ restricts to a map 
$\Hom_{\CO} (U,V) \to \Hom_{A}(U,V)$. By Higman's criterion 
for symmetric algebras (cf. \cite{Br09}), we have 
$\Tr^A_1(\Hom_\CO(U,V))=$ $\Hom_A^\pr(U,V)$. Denote by 
$$\varphi:  K\Hom_{\CO}(U, V) \times K\Hom_{\CO}(V,U)\to K$$ 
the $K$-linear map sending $(\alpha,\beta)\in$ 
$\Hom_{\CO}(U, V)\times \Hom_{\CO}(V, U)$ to 
$\tr_U(\beta\circ\alpha)$, and  denote by 
$$\varphi_A : K\Hom_{A}(U, V) \times K\Hom_{A}(V, U)\to K$$ 
the map sending $(\alpha,\beta)\in$ 
$\Hom_{A}(U, V)\times \Hom_{A} (V, U)$ to 
$\tr_{KU}(z_A^{-1}\beta\circ\alpha)$, where $\alpha$, $\beta$
are extended linearly to maps between $KU$, $KV$.
The following fact generalises \cite[Prop.~1.1]{Thev88Dual}.

\begin{Proposition} \label{Adjunctions} 
With the notation above, for $\alpha\in$ $\Hom_\CO(U,V)$ and
$\beta\in$ $\Hom_A(V,U)$ we have
$$\varphi_A(\Tr^A_1(\alpha),\beta)=\varphi(\alpha,\beta)\ .$$
Similarly, for $\gamma\in$ $\Hom_A(U,V)$ and $\delta\in$ 
$\Hom_\CO(V,U)$ we have
$$\varphi_A(\gamma,\Tr^A_1(\delta))=\varphi(\gamma,\delta)\ .$$
In particular,  $\varphi_{A} $ is non-degenerate.
\end{Proposition}

\begin{proof}      
We regard $\Hom_\CO(U,V)$ and $\Hom_\CO(V,U)$ as $A$-$A$-bimodules
in the canonical way.
If $\mu \in \Hom_{\CO}(U,V)$ and $\beta\in\Hom_{\CO}(V, U)$, 
then for any $a\in A$, we have $\beta \circ a\mu  =$
$\beta a \circ \mu$. If $\epsilon\in\End_{\CO}(U)$ and $a \in A$, 
then $\tr_U(\epsilon a) =$ $\tr_U(a \epsilon)$. Thus we have
\begin{align*} \varphi_A(\Tr^A_1(\alpha),\beta)
& =
\tr_{KU}(z_A^{-1}\sum_{x\in X}\beta\circ x\alpha x^\vee)=
\tr_{KU}(z_A^{-1}\sum_{x\in X}x^\vee\beta\circ x\alpha) \\
& =
\tr_{KU}(z_A^{-1}\sum_{x\in X}x^\vee \beta x\circ \alpha)=
\tr_{KU}(z_A^{-1}\sum_{x\in X}x^\vee x\beta \circ \alpha)\\
&=
\tr_{KU}(z_A^{-1} z_A\beta \circ \alpha)=
\varphi(\alpha,\beta)\ .
\end{align*}
This shows the first equality, and the proof of the second is
analogous. Clearly $\varphi$ is non degenerate, and hence so is 
$\varphi_A$.
\end{proof}

\begin{proof}[{Proof of Theorem \ref{TatedualityO}}]
For $E$ an $\CO$-submodule of $\Hom_{KA}(KU,KV)$ denote
by $E^\perp$ the $\CO$-submodule in $\Hom_{KA}(KV,KU)$
consisting of all $\beta\in$ $\Hom_{KA}(KV,KU)$ such that
$\varphi_A(\epsilon,\beta)\in$ $\CO$ for all $\epsilon\in$ $E$.
By the previous proposition, $\varphi_A$ is
non degenerate, and hence if $E$ is a lattice in 
$\Hom_{KA}(KU,KV)$, then $E^\perp$ is a lattice in 
$\Hom_{KA}(KV,KU)$, and we have $(E^\perp)^\perp=$ $E$.
We need to show that
$(\Hom_A^\pr(U,V))^\perp=$ $\Hom_A(V,U)$.
Let $\beta\in$ $\Hom_{KA}(KU,KV)$. We have
$\beta\in$ $(\Hom^\pr_A(U,V))^\perp$ if and only if
$\varphi_A(\Tr^A_1(\alpha), \beta)\in $ $\CO$ for all $\alpha\in$
$\Hom_\CO(U,V)$. By Proposition \ref{Adjunctions}, this is
equivalent to $\tr_{KU}(\beta\circ\alpha)\in$ $\CO$ for all
$\alpha\in$ $\Hom_\CO(U,V)$. This, in turn, is the case if and
only if $\beta$ belongs to the subspace $\Hom_A(U,V)$ of
$\Hom_{KA}(KU,KV)$. (To see this, choose a basis of $U$, a
basis of $V$, and let $\alpha$ range over the maps sending
exactly one basis element in $U$ to a basis element in $V$ and
all other basis elements of $U$ to $0$).
\end{proof}
 
\begin{proof}[Proof of Corollary \ref{TateDuality}]
We have $z_A=$ $\pi^n 1_A$. The non degenerate pairing
$$\Hombar_A(U,V)\times\Hombar_A(V,U)\to K/\CO$$ 
from Theorem \ref{TatedualityO} has image contained in the submodule
$\pi^{-n}\CO/\CO$ of $K/\CO$. Multiplication by $\pi^n$
yields an isomorphism $\pi^{-n}\CO/\CO\cong$ $\CO/\pi^n\CO$.
Thus Corollary \ref{TateDuality} follows from Theorem
\ref{TatedualityO}.
\end{proof}

In order to prove Theorem \ref{GenKnorr}, we need the
following generalisation of \cite[Prop.~4.2]{CarJo}.

\begin{Proposition} \label{ConstantValue} 
Let $A$ be a symmetric $\CO$-algebra with symmetrising form $s$
such that $K\tenO A$ is separable. Set $z=z_A$.
Let $U$ be an $A$-lattice and let $a$ be the smallest non-negative 
integer such that $\pi^a $ annihilates $\Endbar_A(U) $.  Then 
$$\pi^a\tr_{KU}(z^{-1}\End_A(U)) = \CO\ .$$
\end{Proposition}  

\begin{proof}
Let $\alpha\in$ $\End_A(U)$. By the assumptions we have
$\pi^a\alpha\in$ $\End^\pr_A(U)$. Theorem \ref{TatedualityO},
applied with $U=V$ and $\beta=\Id_U$ implies that
$\pi^a\tr_{KU}(z^{-1}\alpha)\in$ $\CO$. Thus 
$\pi^a\tr_{KU}(z^{-1}\End_A(U))\subseteq$ $\CO$. For the reverse
inclusion, consider first the case that $U$ is non-projective.
Then $a\geq$ $1$, and $\pi^{a-1}\Id_U$ is not contained in
$\End_A^\pr(U)$; equivalently, its image in $\Endbar_A(U)$ is
nonzero. Again by Theorem \ref{TatedualityO}, there exists
$\alpha\in$ $\End_A(U)$ such that $\pi^{a-1}\tr_U(z^{-1}\alpha)
\notin$ $\CO$. Thus $\pi^a\tr_{KU}(z^{-1}\End_A(U))$ is not
contained in $\pi\CO$, whence the equality in this case.
Suppose $U$ is projective, so $a=0 $. Let $\alpha \in$
$\End_{\CO}(U)$ be such that $(\tr_U(\alpha)) =1$ and set $\beta =$
$\Tr_{1}^A (\alpha) \in$ $\End_A(U)$. By 
Proposition~\ref{Adjunctions}, we have 
$\tr_U(z^{-1}\beta)=$ $\varphi_A(\Tr^A_1(\alpha),\Id_U)=$
$\varphi(\alpha,\Id_U)=$ $\tr_U(\alpha)=1$. The result follows.
\end{proof}

\begin{proof}[{Proof of Theorem \ref{GenKnorr}}]
Let $a$ be the smallest positive integer such that $\pi^a$ annihilates 
$\Endbar_A(U)$. The algebra $\Endbar_A(U)$ is local, as $U$ is 
indecomposable non projective. The duality in Theorem 
\ref{TatedualityO} implies that $\soc(\Endbar_A(U))$ is simple.

Suppose that (i) holds. We show first that $U$ is absolutely
indecomposable. The inequality in (i) applied to the endomorphism
$\alpha$ given by multiplication with $z$ shows that 
$$\nu(\rk_\CO(U))=\nu(\tr_U(\Id_U))\geq \nu(\tr_{KU}(z^{-1}\Id_U))\ ,$$ 
so in particular, $\tr_{KU}(z^{-1}\Id_U)$ is nonzero. 
The inequality in (i) applied to an arbitrary $\alpha\in$ $\End_A(U)$
implies that the scalar $\tau$ defined by
$$\tau = \tr_{KU}(z^{-1}\alpha) \tr_{KU}(z^{-1}\Id_U)\inv$$
belongs to $\CO$. A trivial verification shows that
$$\tr_{KU}(z^{-1}(\alpha-\tau\Id_U)) = 0\ .$$
Thus condition (i) implies that $\alpha-\tau\Id_U$ is not
an automorphism, hence in $J(\End_A(U))$. It follows that
$\End_A(U)=$ $\CO\cdot\Id_U + J(\End_A(U))$, and hence $U$
is absolutely indecomposable.

We show next that $U$ has the stable exponent property.
Since the socle of $\Endbar_A(U)$ is
simple, we have $\soc(\Endbar_A(U))\subseteq$ $\pi^{a-1}\Endbar_A(U)$,
and it suffices therefore to show that $\pi^{a-1}\Endbar_A(U)$
is a semisimple $\Endbar_A(U)$-module. That is, it suffices
to show that $\pi^{a-1}\Endbar_A(U)$ is annihilated by
$J(\Endbar_A(U))$. Let $\alpha\in$ $J(\End_A(U))$. The assumptions 
in (i) together with Proposition \ref{ConstantValue} imply that
$\pi^{a}\tr_{KU}(z^{-1}\alpha)\in$ $\pi\CO$, hence
$\pi^{a-1}\tr_{KU}(z^{-1}\alpha)\in$ $\CO$. By Theorem
\ref{TatedualityO}, this is equivalent to 
$\pi^{a-1}\alpha\in$ $\End_A^\pr(U)$, or equivalently, to
$\pi^{a-1}\ualpha=0$. This shows that (i) implies (ii).

Suppose conversely that (ii) holds.
In particular, the socle of $\Endbar_A(U)$ is simple and equal to
$\pi^{a-1}\Endbar_A(U)$. Let $\alpha\in$ $J(\End_A(U))$. The
image $\ualpha$ in $\Endbar_A(U)$ is contained in $J(\Endbar_A(U))$,
and hence $\ualpha$ annihilates $\pi^{a-1}\Endbar_A(U)$. 
Thus $\pi^{a-1}\ualpha=0$. Theorem \ref{TatedualityO} implies that
$\pi^{a-1}\tr_{KU}(z^{-1}\alpha)\in$ $\CO$, hence
$\pi^{a}\tr_{KU}(z^{-1}\alpha)\in$ $\pi\CO$.

By Proposition \ref{ConstantValue}, there exists $\alpha\in$
$\End_A(U)$ such that $\pi^{a}\tr_{KU}(z^{-1}\alpha)= $ $1$.
By the previous argument, this forces $\alpha\notin$
$J(\End_A(U))$. Since $U$ is absolutely indecomposable, it follows
that $\End_A(U)$ is split local, and hence we have $\alpha=$ 
$\lambda\Id_U + \rho$ for some $\lambda\in$ $\CO^\times$ and some 
$\rho\in$ $J(\End_A(U))$.
Since $\pi^a\tr_{KU}(z^{-1}\rho)\in$ $\pi\CO$, it follows
that $\pi^a\tr_{KU}(z^{-1}\lambda\Id_U)\in$ $\CO^\times$.
Then in fact  $\pi^a\tr_{KU}(z^{-1}\lambda\Id_U)\in$ 
$\CO^\times$ for any $\lambda\in$ $\CO^\times$, and hence
$\pi^a \tr_{KU}(z^{-1}\alpha)\in$ $\CO^\times$ for any
automorphism $\alpha$ of $U$. This shows that (ii) implies (i).
\end{proof}

\begin{proof}[{Proof of Theorem \ref{KnoerrExponent}}]   
Let $n$ be the positive integer such that $z_A=$ $\pi^n 1_A$, for
some choice of a symmetrising form. Condition (i) in 
Theorem \ref{GenKnorr} is then equivalent to stating that
$U$ is a Kn\"orr lattice. Thus Theorem \ref{KnoerrExponent}
follows from Theorem \ref{GenKnorr}.
\end{proof} 

%%%%%%%%%%%%%%%%%%%%%%%%%%%%%%%%%%%%%%%%%%%%%%%%%%%%%%%%%%%%%
\section{Characterisations  of the projective scalar
property}
\label{sec_el_char}
 
Throughout this section, $A$ will denote an $\CO$-order such that 
$K \tenO A $ is separable. We identify $A$  with  its canonical image 
in $KA=$  $K\tenO A$. 
Denote by $\Irr_K(A)$ the set of the characters of the simple
$KA$-modules. For $\chi\in \Irr_K(A)$ denote by $e(\chi)$ the
unique primitive idempotent in $Z(KA)$ satisfying $\chi(e(\chi))\neq$ 
$0$. We will use this notation for other orders as well.

\begin{proof}[{Proof of Theorem \ref{scalarA}}]  
Suppose that  $K\tenO A $ is split semisimple.
Proposition \ref{ConstantValue} shows that (i) implies (iii).

By the assumptions, $KAe(\chi)$ is a
matrix algebra over $K$ of dimension $\chi(1)^2$. In particular,
$KAe(\chi)$ is symmetric with symmetrising form $\chi$, and we have
$Z(KA)=$ $\prod_{\chi\in \Irr_K(A)} Ke(\chi)$. Fix a symmetrising form
$s$ of $A$. Then $s$ extends to a symmetrising form of $KA$,
still denoted $s$, and we have
$$ s = \sum_{\chi\in \Irr_K(A)} \sigma_\chi \cdot \chi$$
for some $\sigma_\chi\in$ $K$. The relative projective element
of the matrix algebra $KAe(\chi)$ with respect to $\chi$ is
$\chi(1)\cdot e(\chi)$, and hence the relative projective element of 
$A$ with respect to $s$ is
$$z_A = \sum_{\chi\in \Irr_K(A)} 
\sigma_\chi^{-1}\cdot \chi(1)\cdot  e(\chi)\ .$$
Suppose that (ii) holds; that is, we may assume that $s$ 
satisfies
$$s = \pi^{-n} \rho = 
\sum_{\chi\in \Irr_K(A)} \pi^{-n}\cdot \chi(1) \cdot \chi\ .$$
In that case, a trivial calculation shows that the associated 
relative projective element is, by the previous formula, equal to 
$\pi^n \cdot 1_A$. Thus (ii) implies (i). 

Suppose finally that (iii) holds. We need to show that (ii) holds.
For $U$ an $A$-lattice, write as before $KU=$ $K\tenO U$, and denote 
by $a(U)$ the smallest non-negative integer such that $\pi^{a(U)}$
annihilates $\Endbar_A(U)$. By the assumptions in (iii) and by
Proposition \ref{ConstantValue}, there is
a non-negative integer $n$ such that
$$\pi^n\cdot \tr_{KU}(z^{-1}_A\cdot \End_A(U)) = \tr_U(\End_A(U))=
\pi^{n-a(U)}\CO$$
for any $A$-lattice $U$. We apply this first to $U=A$.
Since $A$ is projective as a left $A$-module, we have $a(A)=0$,
and hence
$$\pi^n\cdot \tr_{KA}(z^{-1}_A\cdot \End_A(A)) = 
\tr_A(\End_A(A)) = \pi^n\CO$$
Any $A$-endomorphism is given by right multiplication with
an element $a$ in $A$. By elementary linear algebra, the trace
of this endomorphism is equal to the trace of the linear
endomorphism given by left multiplication with $a$, and hence
this trace is equal to $\rho(a)$.
Thus $\tr_A(\End_A(A))=$ $\rho(A)=$ $\pi^n\CO$, which implies that
$\pi^{-n}\rho$ sends $A$ to $\CO$. Thus we have
$$\pi^{-n}\rho = s_w$$
for some $w\in$ $Z(A)$. In order to show that $\pi^{-n}\rho$ is a
symmetrising form on $A$ we need to show that $w\in$ $Z(A)^\times$.
Writing $w = \sum_{\chi\in \Irr_K(A)} \omega_\chi e(\chi)$ with
coefficients $\omega_\chi\in$ $\CO$, we need to show that
$\omega_\chi\in$ $\CO^\times$. A trivial calculation shows that
$$s_w = \sum_{\chi\in \Irr_K(A)} \sigma_\chi\omega_\chi \chi$$
Comparing coefficients with $\pi^{-n}\rho$ yields therefore
$$\sigma_\chi\omega_\chi = \pi^{-n}\chi(1)$$
for all $\chi\in \Irr_K(A)$, and hence
$$\nu(\sigma_\chi\omega_\chi) = \nu(\pi^{-n}\chi(1))\ .$$
Let $\chi\in \Irr_K(A)$, and let $V$ be an $A$-lattice such that
$KV=$ $K\tenO V$ has character $\chi$. Using that
$\End_A(V)=$ $\CO\cdot\Id_V$, we get from the above that
$$\nu(\pi^n\cdot \tr_{KV}(z_A^{-1})) = \nu(\Id_V) = \nu(\chi(1))\ .$$
By the above formula for $z_A$, we have $z_A^{-1}=$ 
$\sum_{\chi\in \Irr_K(A)} \sigma_\chi\cdot  \chi(1)^{-1}\cdot e(\chi)$,
and hence $\tr_{KV}(z_A^{-1}) = \sigma_\chi$. Thus
$$\nu(\pi^n\sigma_\chi) = \nu(\chi(1))$$
Combining the previous statements yields
$$\nu(\sigma_\chi\omega_\chi)= \nu(\pi^{-n}\chi(1)) =
\nu(\sigma_\chi)$$
and hence $\omega_\chi$ is invertible in $\CO$. This
shows that (iii) implies (ii).
The last statement in Theorem \ref{scalarA} on the integer
$n$ is obvious from the proofs of the implications.
\end{proof}

\begin{Remark}
The coefficients $\sigma_\chi^{-1}$ in the above proof are called 
{\it Schur elements} in \cite[\S 7.2]{GP}. 
\end{Remark}
 
Next, we prove Theorem~\ref{MoritascalarA}. As in the theorem, let 
$\Irr_k(A)$ denote an indexing set for the isomorphism classes of 
simple  $k\tenO A $-modules, and for $ \chi  \in  \Irr_K(A)$ and 
$ \varphi \in \Irr_k(A) $ denote by $ d_{\chi, \varphi} $ the   
multiplicity  of $S$ as a composition  factor of $k\tenO V$, where  
$V$ is an $A$-lattice  such that $ K \tenO V$ has character 
$\chi $, and $S$ is a simple $k\tenO A$-module with character  
$\varphi$.
We adopt the analogous notation for other orders.

\begin{Lemma} \label{deriveddeco}  
Let $A'$  be an $\CO$-order   which is derived equivalent to  $A$. 
Then $|\Irr_k(A)| = |\Irr_k(A')|$, and $|\Irr_K(A)| =|\Irr_K(A')| $.  
Further,  there exists a bijection   $\chi \to \chi' $    from  
$\Irr_K(A)$ to $\Irr_K(A')$,   signs 
$\epsilon_{\chi } \in  \{\pm 1 \}$, $\chi \in   \Irr_K(A)$, 
and integers  $u_{\varphi, \psi} $, $ \varphi \in \Irr_k(A)$, 
$ \psi \in \Irr_k(A')$ such that

\begin{enumerate}
\item[\rm(i)]    
For $\chi \in \Irr_K(A),  \psi \in \Irr_k(A') $,    
$d_{\chi', \psi}  =  
\epsilon_{\chi}\sum _{\varphi\in \Irr_k(A) } 
d_{\chi, \varphi}u_{\varphi, \psi} $.

\item[\rm(ii)]    
The form  $s= \sum _{\chi \in \Irr_K(A)}  \sigma_{\chi} \chi $, 
$\sigma_{\chi} \in K$ is a symmetrising form of $A$ if and only if 
the form  $s'=\sum _{\chi \in \Irr_K(A)}
\epsilon_{\chi} \sigma_{\chi} \chi' $ is a symmetrising form of $A'$. 
\end{enumerate}

If $A$ and $A'$ are Morita equivalent, then in addition there is 
a bijection $ \varphi \to \varphi'$ from $\Irr_k(A)$ to 
$\Irr_k(A')$ such that  $d_{\chi', \varphi'} =  d_{\chi, \varphi} $   
and $\epsilon_{\chi}=1 $  for all $\chi \in \Irr_K(A) $, 
$ \varphi\in \Irr_k(A)$.
 
\end{Lemma}

\begin{proof}  
The first  statement follows from  
\cite[Theorem 6.8.8]{Zimmermannbook}. The transfer  of symmetrising 
forms  as in (ii)  is proved in \cite[Theorem 4.7]{Eisele12}.
\end{proof}

\begin{proof}[{Proof of Theorem~\ref{MoritascalarA}}]      
Suppose that the $\CO$-order $A'$ is Morita equivalent to $A$ and 
let $\chi \to \chi' $, and $\varphi \to \varphi'$ be the bijection  
of  Lemma~\ref{deriveddeco}.  Denoting by $n_{\varphi}$ the 
$k$-dimension of the simple $A'$-module labelled by $\varphi'$ 
($\varphi \in \Irr_k(A)$), we have that $\chi' (1) =$ 
$\sum_{\varphi \in \Irr_k(A)} n_\varphi\cdot d_{\chi, \varphi}$ for 
all $\chi \in \Irr_K(A) $. The equivalence between (i) and (iii) is 
now immediate from Lemma~\ref{deriveddeco}  and  the equivalence  
between (i) and (ii)  of Theorem~\ref{scalarA}.
We now prove that (iv) implies (iii). Let $n$ and $m_{\varphi}$, 
$\varphi \in \Irr_k(A)$ be integers such that  
$\pi^{-n} \sum_{ \chi  \in \Irr_K(A)} a_{\chi} \chi$ is a   
symmetrising form of $A$, where $a_{\chi} =$ 
$\sum _{\varphi \in \Irr_k(A)} m_\varphi d_{\chi, \varphi} $, 
$ \chi \in \Irr_K(A)$.  
Let $X$ be an $\CO$-basis of $A$. Choose a positive integer $t$  such 
that
$ \pi^{-n} \cdot p^t \cdot d_{\chi, \varphi} \chi(x) \in \pi \CO$  
and $m'_\varphi:=m_\varphi +p^t > 0$ for all $\chi \in \Irr_K(A)$, 
$ \varphi \in \Irr_k(A)$ and $x\in X$. Set $s'=$
$\pi^{-n} \sum_{\chi \in \Irr_K(A)} a_{\chi}'\cdot \chi$, where 
$ a_{\chi}'=$ 
$\sum_{\varphi \in \Irr_k(A)} m_\varphi'\cdot  d_{\chi, \varphi} $, 
$\chi \in  \Irr_K(A) $. Then for all $a\in A$,  
$s'(a)-s(a) \in \pi \CO $.  Hence by considering the 
determinant of the Gram matrices of the bilinear forms associated to  
$s$ and $s'$, it follows that $s'$ is also a symmetrising form of 
$A$. This proves that (iii) holds.
Since (i) clearly implies (ii) and (iii) implies (iv), in order to 
complete the proof, it suffices to show that (ii) implies (iv).  
Suppose that $A'$ has the scalar projective property and  that $A'$ 
and $A$ are derived equivalent. Then  (iii) holds for $A'$, say for 
the integers $m_{\psi}$, $\psi \in \Irr_k(A')$. Then by 
Lemma ~\ref{deriveddeco}, we have that (iv) holds for $A$ with the  
integers  $n_\varphi = \sum_{\psi \in \Irr_k(A')} m_{\psi}u_{\varphi, \psi}$, 
$\varphi \in  \Irr_k(A)$.
\end{proof}

For the rest of this section we will expand on the question of the 
extent to which the characterisations of the projective scalar 
property up to Morita equivalence given in Theorem \ref{MoritascalarA} 
(iii) and (iv) are constructive. The point here is that the set of 
symmetrising forms for an order $A$ is actually a $Z(A)^\times$-orbit, 
and $Z(A)$ is an $\CO$-order for a (potentially) quite large ring 
$\CO$. But in fact, as we will see, the criterion can be reduced to 
linear algebra over $\Q$. 

The following proposition shows that the projective scalar property is 
essentially independent of the choice of the ring $\CO$. This is 
particularly interesting to note since we often make the assumption 
that $K$ is a splitting field.

\begin{Proposition}
Let $A$ be an $\CO$-order and let $\mathcal E\supseteq \CO$ be a 
discrete valuation ring containing $\CO$ such that 
$J(\mathcal E) \cap \CO = J(\CO)$. Then an $\CO$-order $A$ has the 
projective scalar property if and only if the $\mathcal E$-order 
$\mathcal E \tenO A$ has the projective scalar property.
\end{Proposition}

\begin{proof}
By the characterisation in Theorem \ref{scalarA}, $A$ having the 
projective scalar property is equivalent to some multiple of the 
regular trace being a symmetrising form for $A$. But the regular trace 
on $A$ and the regular trace of $\mathcal E \tenO A$ have the same 
Gram-matrix (when the same basis is chosen for both of them), and 
invertibility of a multiple of said Gram-matrix over $\CO$ is 
equivalent to invertibility over $\mathcal E$, provided of course 
that we multiplied by an element of $\CO$.
	
So the only thing that still requires proof is that if $\tau$ is a 
generator of $J(\mathcal E)$, then the integer $m$ such that 
$\tau^{-m}\cdot  \rho$ is a symmetrising form for $\mathcal E \tenO A$ 
satisfies $\tau^{m} \mathcal E = \pi^{n} \mathcal E$ for some 
$n \in \Z_{\geq 0}$ (since this means that $\pi^{-n}\cdot \rho$ is a 
symmetrising form for $A$). But by Theorem \ref{scalarA} we have
$\tau^m \mathcal E =$ 
$\tr_{\mathcal E \tenO A} (\End_{\mathcal E \tenO A} 
(\mathcal E \tenO A)) = \mathcal E \tenO 
\tr_A(\End_A(A))$, 
and $\tr_A(\End_A(A))$ is certainly of the form $\pi^n\CO$ for some $n$.
\end{proof}

\newcommand{\Zrat}{Z^{\rm rat}}

\begin{Definition}
Let $A$ be an $\CO$-algebra which is free of finite rank as
an $\CO$-module such that $KA$ is split semisimple.
Fix an isomorphism 

$$\varphi:\  Z(K A) \stackrel{\sim}{\longrightarrow}  
K \times \ldots \times K $$
  
We define the \emph{rational centre} $\Zrat(KA)$ of 
$KA$ to be the $\Q$-algebra 
$$ \varphi^{-1}(\Q\times \ldots\times \Q) $$ 
We define the {\it rational centre of $A$}, denoted by $\Zrat(A)$, as 
the intersection of $A$ with $\Zrat(KA)$.

We say that $A$ is \emph{rationally symmetric} if there is an
element 
$$\tilde \sigma = 
\sum_{\chi \in \Irr_K(A)}\tilde \sigma_\chi e_\chi \in \Zrat(A)$$
and an $n\in \Z$ such that 
$$ \pi^{-n}\cdot \sum_{\chi \in \Irr_K(A)} 
\tilde \sigma_\chi\cdot \chi $$ 
is a symmetrising form for $A$.
\end{Definition}

We should note that $\sigma_\chi = \pi^{-n}\cdot \tilde \sigma_\chi$ 
with $\sigma_\chi$ defined as earlier. Therefore rational symmetry is 
not the same as asking that the $\sigma_\chi$ be rational. Not even 
the projective scalar property implies rationality of the 
$\sigma_\chi$. 

The rational centre of $A$ is a $\Z_{(p)}$-order, and the projective 
scalar property implies rational symmetry. 
We should remark that, if $\CO$ is ramified over $\Z_p$, then rational 
symmetry is not necessarily preserved under direct sums. Neither is 
the projective scalar property, or even the property of being 
Morita-equivalent to an order which satisfies the projective scalar 
property. This is due to the possibility that the rational 
symmetrising forms involve different powers of $\pi$, whose quotient 
may have a non-integral $p$-valuation (using the convention $\nu(p)=1$). 

\begin{Remark}
An element $\tilde\sigma$ (together with an $n \in \Z$) as above  
and  the central projective element $z_A$ are related by the formula
$$ z_A = \pi^n\cdot \tilde \sigma^{-1}\cdot 
\sum_{\chi\in \Irr_K(A)} \chi(1) \cdot e_\chi $$ 
In particular, $\tilde \sigma$ can be chosen in $\Zrat(A)$ if and only 
if $z_A \in K^\times\cdot \Zrat(A)$. Now we can reinterpret the 
projective scalar property and rational symmetry in the following 
way: we consider the orbit $Z(A)^\times \cdot z_A$. If it intersects 
non-trivially with $K^\times \cdot \Zrat(A)$, then $A$ is rationally 
symmetric, and if it intersects non-trivially with 
$K^\times \cdot 1_A$, then $A$ has the projective scalar property.
\end{Remark}

In view of everything we have seen so far, the following is fairly 
straight-forward.

\begin{Proposition}
Assume that $A$ is rationally symmetric, and 
$\tilde \sigma \in \Zrat(A)$ is as before.
Then $A$ has the scalar projective property if and only if

\begin{equation}\label{eqn intersection 1}
     \left\langle  \sum_{\chi \in \Irr_K(A)} 
      d_{\chi, \varphi} \cdot \chi\ \Big|\ 	
      \varphi \in \Irr_k(A) \right\rangle_\Q \cap 
     \left\{ \sum_{\chi\in\Irr_K(A)}
     \tilde \sigma_\chi \cdot \frac{\chi(z)}{\chi(1)} 
     \cdot \chi\ \Big|\ z \in \Zrat(A) \right\} 
\end{equation}
	properly contains
\begin{equation}\label{eqn intersection 2}
     \left\langle  \sum_{\chi \in \Irr_K(A)} 
     d_{\chi, \varphi} \cdot \chi\ \Big|\ \varphi\in\Irr_k(A)  
     \right\rangle_\Q \cap 
     \left\{ \sum_{\chi\in\Irr_K(A)}\tilde \sigma_\chi \cdot 
     \frac{\chi(z)}{\chi(1)} \cdot \chi\ \Big|\ z \in I \right\} 
\end{equation}
for all maximal Ideals $I$ in $\Zrat(A)$. 

\end{Proposition}

Note that the right hand 
side in both \eqref{eqn intersection 1} and \eqref{eqn intersection 2} 
is the intersection of a $\Q$-vector space and a $\Z_{(p)}$-lattice, 
which can be computed by means of linear algebra.

We conclude this section with an example of a symmetric algebra which 
is not rationally symmetric, to show that the two notions are not 
equivalent.

\begin{Example}
Assume that $k$ has characteristic two and $\CO$ is unramified, i. e. 
$\pi=p=2$. Let $x\in \CO^\times$ be an arbitrary unit in $\CO$. 
We consider the order 
$A=\langle \lambda_1,\lambda_2, \lambda_3,\lambda_4\rangle_\CO$ in the 
commutative split-semisimple $K$-algebra 
$K\times K \times K \times K$, where
     \begin{equation}\label{matrix of non-rational example}
       \begin{array}{cccccccc} \lambda_1&=&(&1&1&1&1&)\\
       \lambda_2&=&(&0&2&0&2x&) \\ \lambda_3&=&(&0&0&2&2x&)\\
       \lambda_4&=&(&0&0&0&4x&)
       \end{array}
     \end{equation}
	We claim that the map
       \begin{equation}\label{symm form one}
         s:\ K\times K \times K \times K \longrightarrow 
         K: \ (a_1,a_2,a_3,a_4) \mapsto 
         \frac{2-x^{-1}}{4} a_1 + \frac{1}{4} a_2 + \frac{1}{4} a_3 + 
         \frac{x^{-1}}{4} a_4
       \end{equation}
defines a symmetrising form for $A$. The Gram-matrix of $s$ with 
respect to the basis $(\lambda_1,\ldots,\lambda_4)$ is
       \begin{equation}
         (s(\lambda_i\cdot \lambda_j))_{i,j} = 
         \left( \begin{array}{cccc} 1&1&1&1 \\ 1 & 1+x & x & 2x  
         \\ 1 & x & 1+x & 2x \\ 1 & 2x & 2x & 4x \end{array} \right)
       \end{equation}
The determinant of this matrix is congruent to $1$ mod $2\CO$, which 
implies that it is invertible over $\CO$, which in turn implies that 
$A$ is a self-dual lattice with respect to $s$. So clearly, $A$ is a 
symmetric $\CO$-order. However, if $x+2\CO \neq 1 + 2\CO$, then $A$ 
is not rationally symmetric. To see this we consider the family of 
forms 
          \begin{equation}
          s_u:\ K\times  K \times K \times  K \longrightarrow 
          K:\ (a_1,a_2,a_3,a_4) \mapsto \frac{1}{4} \cdot 
          \sum_{i=1}^4 u_i\cdot a_i
          \end{equation}
where $u\in (K^\times)^4$. By definition, the order $A$ is rationally 
symmetric if and only if $s_u$ is a symmetrising form for $A$ for some 
$u\in (\mathbb Q^\times)^4$. 
We know that the symmetrising forms for $A$ are exactly the forms 
$s(z\cdot -)$ with $z\in Z(A)^\times$ and $s$ as in 
\eqref{symm form one}. 
The form $s(z\cdot -)$ is equal to $s_{z\cdot v}$ with 
$v=(2-x^{-1},1,1,x^{-1})$. Since $z$ is a unit each $z_i$ lies in 
$\CO^\times$, and so do all $v_i$. So if $A$ is symmetric with respect 
to $s_u$, then each $u_i$ needs to lie in $\CO^\times$.
Moreover, $A$ being symmetric with respect to $s_u$ would necessitate 
$A$ being integral with respect to $s_u$, which in particular would 
require $s_u(\lambda_2)= 2^{-1}\cdot u_2+2^{-1}\cdot u_4\cdot x \in 
\CO$. That is, $-\frac{u_2}{u_4} + 2\CO =x + 2\CO$, which can only 
hold true for rational $u_i$'s if $x+2\CO$ lies in the prime field of 
$k$, which means $x+2\CO = 1+2\CO$ (since we asked that $x$ be a unit, 
the case $x+2\CO = 0+2\CO$ is impossible). 
\end{Example}

%%%%%%%%%%%%%%%%%%%%%%%%%%%%%%%%%%%%%%%%%%%
\section{Heights and Degrees of    Kn\"orr    lattices} \label{heights}

\begin{proof}[{Proof of Proposition \ref{scalardegree}}]
Let $U$ be  a Kn\"{o}rr lattice.  Then,  $\tr_U(\End_A(U)) = 
\rank(U)\CO $.
By Theorem \ref{scalarA} (iii), we have that 
$\tr_U(\End_A(U)) = \pi^{n-a(U)}\CO $.
By Theorem \ref{scalarA} (ii), we have
$$\frac{\rank(A)}{\pi^n} = \frac{\rho(1_A)}{\pi^n} \in \CO\ .$$
It then follows that
$$ \rank(A)\CO \subseteq \pi^n \CO \subseteq \pi^{n-a(U)}\CO = 
\rank(U)\CO\ .$$ 
This proves (i). Now suppose that $U$ is a projective lattice. Then, 
by Theorem \ref{scalarA} (iii) we have that $\tr_U(\End_A(U)) =
\pi^n \CO $. On the other hand, $\rank (U) \CO \subseteq 
\tr_U(\End_A(U)) $. This proves (ii).
\end{proof}
 
The next lemma is needed to prove Theorem \ref{scalarheight}.

\begin{Lemma}\label{lemmascalarheight}  
Let $A$  be a symmetric $\CO$-algebra such that $ K \tenO A$ is split 
semisimple. Assume that $A$ has the projective scalar property and let 
$\pi^n1_A$ be a relative projective element  with respect to   some 
symmetrizing form on $A$.  Let $a_0=\text{max}_V \{a(V)\} $   as  
$V$ ranges over  all $A$-lattices.  There exists  $\chi \in \Irr_K(A)$ 
such that $\chi(1) \CO = \pi^{n-a_0}\CO$. 
\end{Lemma}
 
\begin{proof} 
By Theorem  \ref{scalarA}(ii), $ a_0 \leq n $ and  we have that 
$\chi(1) \CO \subseteq  \pi^{n-a_0}\CO $ for all $\chi \in \Irr_K(A)$.
Let $U$ be an $A$-lattice and $\alpha \in \End_A(U) $ be such that  
$\tr_U(\alpha)\CO = \pi^{n -a_0}\CO$. Let $ f \in \CO[x] $ be the 
characteristic polynomial of $ \alpha $  and let  $g \in \CO[x]$  be 
an irreducible  monic factor  of $f$.
Let  $\bar K $   be an algebraic closure of $K$,   let $\lambda_i$,   
$i \in {\mathcal I}$ be the roots of $g$  in  $\bar K$  and  for 
$i \in  {\mathcal I}  $  let $W_{i}  $ be the  generalised  
$\lambda_i $-eigen space of  $\alpha $   in  $ \bar K \tenO U $.  Set  
$\lambda_g:= \sum_{i\in {\mathcal I}} {\lambda_i} \in \CO $  and  
$ W_g := \oplus_{i\in {\mathcal I}} W_i $.
We have  $\dim (W_i) = \dim (W_j)  =:d_g $ for all $i$, $j \in 
{\mathcal I} $.  So,
$$ \tr_{W_g}(\alpha) = \lambda_g   d_g\ .$$  
Since $\tr_{U}(\alpha)  $ is the sum  of $ \tr_{W_g}(\alpha) $ as $g$ 
runs through the   irreducible factors of $f$  and since $\lambda_g 
\in \CO $, replacing $g$  by some other irreducible factor of $f$
if necessary, we may assume that 
$$\pi^{n-a_0}\CO = \tr_{U}(\alpha) \CO \subseteq d_g \CO\  .$$
Now $\alpha  \in \End_A (U) $,  hence $ U_i $ is a 
$\bar K \tenO A$-submodule of $\bar K \tenO A$. In particular, $d_g $ 
is the dimension of a $\bar K \tenO A$-module. Since $K \tenO A$ is 
split it follows that there  exists some  $\chi \in I$ such that 
$$\pi^{n -a_0}\CO= \tr_{U}(\alpha) \CO \subseteq 
d_g  \CO  \subseteq \chi(1) \CO  \subseteq \pi^{n -a_0}\CO.$$ 
Hence, $ \chi (1) \CO =  \pi^{n -a_0}\CO$  as desired.
\end{proof} 

\begin{proof} [{Proof of Theorem \ref{scalarheight} }]   
The fact that  $U'$  is a Kn\"{o}rr  $A'$-lattice is a consequence of 
Theorem \ref{KnoerrExponent}.
Let $a_0=\text{max}_V \{a(V)\}$ as  $V$ ranges over all $A$-lattices.
Then  $a_0 $ also equals $\max_{V'} \{a(V')\} $ as  $V'$ ranges over   
all $A'$-lattices. Further, $a(U)=a(U') $.
Let $\pi^e \CO = p\CO$   and let  $\pi^n1_A$ be  a  relative 
projective element of $A$. For any $A$-lattice $V$, 
$ \rank (V) \CO \subseteq \tr_V(\End_A(V) )$, hence 
by Lemma \ref{lemmascalarheight}  and Theorem \ref{scalarA} 
$$p^{\frac{ n-a_0}{e}} = \text{min}_V \{\rank (V )_p\}$$    
as  $V$ ranges over all $A$-lattices. 
Since $U$ is a Kn\"{o}rr $A$-lattice  and using again  Theorem 
\ref{scalarA}, it follows that 
$$\pi^{n-a(U)}\CO  = p^{\frac{ n-a_0}{e} +h(U)}\CO$$  
and hence   
$$h(U)= \frac{a_0 -a(U)}e\ . $$   
Applying the same argument to $A'$ and $U'$ gives the desired result.
\end{proof}

\begin{proof}[{Proof of Proposition \ref{knorrproj}}]
Let $u$ be an element of $U \setminus \pi U.$ Let 
$\varphi: U \rightarrow U$ be an $\CO$-linear projection onto 
$\CO u,$ and let   $\Tr_1^A(\varphi)$  be the corresponding 
$A$-endomorphism of $U$. A  calculation similar to that in Proposition 
\ref{Adjunctions} and using the assumption that $\pi^n 1_A$ is a 
relative projective element as in Theorem \ref{scalarA}, shows that
$$\tr_U(\Tr_1^A(\varphi)) = \pi^n\ . $$
Now because $U$ is projective, we have $a(U)=0$. It follows that
$$\tr_U(\End_A(U)) = \pi^n\CO\ .$$
Because $U$ is a Kn\"{o}rr lattice, we can conclude that 
$\Tr_1^A(\varphi)$ is an invertible element of $\End_A(U)$. In 
particular, it is surjective. However, the image of $\Tr_1^A(\varphi)$ 
is contained in the $A$-lattice $Au.$ We thus have $Au = U$. The 
result follows because $u$ was an arbitrary element of 
$U\setminus \pi U$.
\end{proof}

%%%%%%%%%%%%%%%%%%%%%%%%%%%%%%%%%%%%%%%%%%%%%%%%%%%%%%%%%%%%%%%%%%%%%
\section{Examples} \label{examples}

\begin{Example}  \label{basicexample}    
If $A= \mathrm{Mat}_n(\CO)$ for some positive  integer $n$ or if 
$A=\CO G$ for some finite group  $G$, then $A$  has the scalar 
projective property   (see \cite[Examples and Remarks after Proposition
 3.3] {Br09}).  If  an  $\CO$-algebra $A$  has the projective scalar  
property, and if $B$ is a direct factor of $A$, then $B$   has  the  
projective scalar property. This is immediate from the fact that the 
relative projective element  with respect to a symmetrising form on 
$A$ is independent of the choice of an $\CO$-basis.  If 
$\CO$-algebras $A$ and $B$  have the projective  scalar  property, 
then so   does $A\tenO B$.
However, the  projective  scalar  property is not preserved under 
taking direct products, whilst the  property of being symmetric is. 
For instance  if $p=2$, then  by Proposition~\ref{scalardegree},  
$\CO \times \mathrm {Mat}_2(\CO )$ does not have the  projective 
scalar property. Further, $\CO \times \mathrm {Mat}_2(\CO )$ is Morita 
equivalent to $\CO \times \CO$ from which we see that  the scalar 
projective property is not invariant under Morita equivalence.
\end{Example}

\begin{Example} \label{sourcealgebraexample} 
Source algebras of blocks of finite groups have the projective scalar 
property. More precisely, if $A$ is a source algebra of a block of a 
finite group algebra with defect group $P$, and  $k$ is a splitting 
field for   the underlying finite group and its subgroups, then there 
is a symmetrising form on $A$ such that the relative projective 
element of $A$ is equal to $|P|\cdot 1$. To see this, let $G$ be a 
finite group, $B$ a block algebra of $\OG$, $P$ a defect group of $B$, 
and $i$ a source idempotent of $B$; that is, $i$ is a primitive 
idempotent in $B^P$ satisfying $\Br_P(i)\neq$ $0$, where 
$\Br_P : (\OG)^P\to$ $kC_G(P)$ is the Brauer homomorphism.   Assume 
that $k$ is a splitting field for $G$ and all of its subgroups.
The source  algebra $A=i\OG i$ is again symmetric, and any 
symmetrising form on $\OG$ restricts to a symmetrising form on $A$.
Denote by $s : \OG\to$ $\CO$ the canonical symmetrising form, sending 
$1_G$ to $1_\CO$ and $x\in$ $G\setminus \{1_G\}$ to zero. With
respect to this form, the relative trace $\Tr^{\OG}_1$ on $\OG$ is
equal to the relative trace map $\Tr^G_1$, sending $a\in$ $\OG$ to
$\sum_{x\in G} xax^{-1}$. The relative trace map $\Tr^A_1$ with respect
to the symmetrising form $s$ restricted to $A$ satisfies $\Tr^A_1(a)=$ 
$\Tr^G_1(a)i$. In particular, we have $\Tr^A_1(i)=$ $\Tr^G_1(i)i$.
As a consequence of \cite{PiPu} or \cite[9.3]{Thev88Dual},  
the element $u=$ $\Tr^G_P(i)$ is invertible in $Z(B)$. 
Moreover, we have $\Tr^G_1(i)=$ $|P|\Tr^G_P(i)=$ $|P|u$.
Denote by $t$ the symmetrising form given by $t(a)=$ $s(ua)$.
The relative trace map on $A$ with respect to the form $t$ sends the
unit element $i$ of $A$ to $|P| u u^{-1} i = |P| i$ as required.
\end{Example}

\begin{Example}  \label{hopfexample}    
If  $A$ is a Hopf algebra  over $\CO$ such that $K\tenO  A$ is  
semisimple, then  $A$ has the projective scalar property. This is 
well-known to Hopf algebra experts-we just sketch the trail of ideas.  
By  \cite[Theorem 3.3]{LaRa88I} and \cite[Theorem 4]{LaRa88II}, 
the antipode of $K\tenO A $ and of $K\tenO A^*= (K\tenO A)^*$ has 
order $2$. Hence the same is true for  the antipode  of $A$ and $A^*$.
By the main theorem of \cite{LaSw},  $A$ 
has a non-singular left integral say  $\lambda $. Then $\lambda$ 
is also a non singular left integral for $K\otimes A$. Hence by   
\cite[Props.~3 and ~4]{LaSw} $\epsilon(\lambda) \ne 0$ and    
$A$ is unimodular. Since the antipode of $A^*$ also has order $2$,   
by the second corollary  to Proposition 8 of 
\cite{LaSw}, applied with  the roles of $A$ and $A^*$ reversed, we 
have that if $\Lambda \in A^* $ is a non-singular integral  
($\Lambda $ exists by the main theorem of \cite{LaSw} applied to   
$A^*$), then  $\Lambda$ is a symmetrising form on $A$. Further, 
by \cite[Section~5.3]{Lore}, the  corresponding projective  element is 
a scalar. 
\end{Example}

\begin{Example} \label{rank2example}
This example shows that very few local commutative symmetric
$\CO$-algebras of $\CO$-rank $2$ have the projective scalar property.
Let $A$ be an indecomposable $\CO$-algebra such that 
$K\tenO A =$ $K\times K$; in particular, $A$ is commutative. Then 
there is a unique positive integer $m$ such that 
$A=$ $\{(\alpha,\beta)\in \CO\times \CO\ |\ 
\beta-\alpha\in \pi^m\CO\}=$ 
$\{(\alpha,\alpha+\beta)\ |\ \alpha\in\CO,\ \beta\in \pi^m\CO\}$. 
The algebra $A$ is local commutative and symmetric, with symmetrising
form $s$ sending $(\alpha,\alpha+\beta)\in$ $A$ to $\pi^{-m}\beta$.
We are going to show that $A$ has the projective scalar property
if and only if $p=2$ and $2\in$ $\pi^m\CO$. 

The $\CO$-basis $X=$ $\{(1,1), (0,\pi^m)\}$ of $A$ has, with respect
to $s$, the dual basis $\{(-\pi^m,0), (1,1)\}$. Thus the relative
projective element with respect to the symmetrising form $s$ is
$z_A =$ $(-\pi^m,\pi^m)$. We have $A^\times=$ 
$\{(\alpha,\alpha+\pi^m\gamma)\ | \ 
\alpha\in\CO^\times, \gamma\in\CO\}$. 
Thus the $A^\times$-orbit of $z_A$ is 
$\{-\pi^m\alpha,\pi^m\alpha+\pi^{2m}\gamma\ |\ 
\alpha\in\CO^\times, \ \gamma\in\CO\}$.
An element in this set is a scalar if and only if
$\pi^m\gamma=$ $-2\alpha$. For $p$ odd this is impossible as the
right side is invertible in $\CO$ whereas the left side has a positive
valuation of at least $m$. This shows that for $p$ odd, $A$ does not
have the projective scalar property. For $p=2$, 
the algebra $A$ has the scalar property if and only if
$\pi^m$ divides $2$ in $\CO$.   

Note that since $A$ is local, any $\CO$-algebra Morita equivalent to 
$A$ is a matrix algebra over $A$. Hence if $A$ does not have  the 
projective scalar property, then neither does any algebra Morita 
equivalent to $A$.
\end{Example}

\begin{Example} \label{heckeexample}
Let $(W,S)$ be a finite Coxeter group with length function $\ell$ and 
$q\in$ $\CO^\times$. Let $\CH = \CH_q(W,S)$ be 
the associated Iwahori-Hecke algebra over $\CO$ with parameter $q$. 
That is, $\CH$ has an $\CO$-basis $\{T_w\}_{w\in W}$, with 
multiplication given by
$T_wT_y=$ $T_{wy}$ if $w$, $y\in$ $W$ such that $\ell(wy)=$ 
$\ell(w)+\ell(y)$, and $(T_s)^2=$ $qT_1 + (1-q)T_s$ for $s\in$ $S$.
By \cite[Proposition 8.1.1]{GP}, the algebra $\CH$ is symmetric, with
a symmetrising form sending $T_1$ to $1$ and $T_w$ to $0$ for
$w\in$ $W\setminus \{1\}$. The dual basis of $\{T_w\}_{w\in W}$
with respect to this form is $\{q^{-\ell(w)}T_{w^{-1}}\}_{w\in W}$, 
and hence the associated relative projective element is 
$$z_\CH = \sum_{w\in W}\ q^{-\ell(w)}T_wT_{w^{-1}}$$
Whether $\CH$ has the projective scalar property seems to be difficult
to read off this expression. 
If $p=2$ and $W=S_2=S=\{1,s\}$, and if $q$ is an odd integer, then
the map sending $T_1$ to $(1,0)$ and $T_s$ to $(1,1-q)$ is an
injective algebra homomorphism from $\CH$ to $\CO\times \CO$. The 
previous example shows that $\CH$ has the scalar property if and only 
if $q\equiv 3\ \mod\ 4$.
\end{Example}

\begin{Example}\label{characterringexample}
Let $G$ be a finite group  and  assume that $\CO$ contains the values 
of all irreducible characters of $G$. Let  
$A=\CO[\Irr(G)] = \CO\otimes_\Z \Z[\Irr(G)]$.  The irreducible 
characters of $G$ form an $\CO$-basis for $A$.  For $K$-valued 
functions $\alpha$ and $\beta$ on $G$, define the usual
\[
[\alpha,\beta] = 
\frac{1}{|G|}\sum_{g \in G} \alpha(g)\beta(g\inv) \in K.
\]
For $\chi \in \Irr(G)$, let $\overline{\chi}$ denote the character of 
the contragredient representation, so $\overline{\chi}(g) = 
\chi(g\inv)$ for all $g\in G$. Finally, let $1_G$ denote the trivial 
character of $G$.

For $\chi$, $\psi \in \Irr(G)$, the identity 
$[\chi\overline{\psi}, 1_G] = [\chi, \psi] = \delta_{\chi, \psi}$ 
implies that the $\CO$-linear function $s: A \rightarrow \CO$
given by  $s(\alpha) = \text{(coefficient of $1_G$ in $\alpha$)} $
is a symmetrizing form on $A$. The same identity makes it clear that 
the basis of $A$ dual to $\Irr(G)$ with respect to $s$ is given by 
$\chi^\vee = \overline{\chi}$. The corresponding relative projective 
element, $z= \sum_{\chi \in \Irr(G)} \chi\overline{\chi}$,
coincides with the function on $G$ sending $g$ to $|C_G(g)|$. Clearly, 
$z$ is a scalar multiple of $1_G$ if and only if $G$ is abelian. In 
this case, we have $z = |G|\cdot 1_G$.
However, to see exactly when $A$ has the projective scalar property, 
it is necessary to consider the action of $A^\times$ on $z$. 
Let $u$ be an invertible element of $A$. Then $u$ is a function from 
$G$ to $\CO^\times$. Assume that $uz = \lambda \cdot 1_G$ for some 
element $\lambda \in \CO$. We must then have
\begin{equation}\label{unitcond}
u(g) = \frac{\lambda}{|C_G(g)|} \in \CO^\times
\end{equation}
for all $g \in G$. Thus, $|C_G(g)|_p$ is independent of $g$. We 
deduce that every element of $G$ must centralise a Sylow $p$-subgroup.
So, let $P$ be a Sylow $p$-subgroup. By Sylow's theorem, every 
element of $G$ is conjugate to an element of $C_G(P)$. A well known 
application of ``Burnside's counting lemma" allows us to conclude that 
$C_G(P) = G$. Thus, $P$ is abelian, and $G \cong P\times H$ for some 
group $H$ of $p'$ order. Conversely, we claim that if
$G=P \times H $, with $P$ an abelian $p$-group and $H$ a 
group of $p'$-order, then $A$ has the projective scalar property.  
All that remains to do is to verify that the function
$u(g) = \frac{1}{|C_G(g)|_{p'}}$ for $g\in G$ 
actually lies in $A$, assuming $G=P\times H$ as above. So let 
$\chi \in \Irr(G)$. We must show that $[\chi, u] \in \CO$. We can 
write $\chi = \theta\otimes \psi$ for irreducible characters $\theta$ 
of $P$ and $\psi$ of $H$. One verifies
$$[\chi,u]= 
   \left\{\begin{array}{cc}
     \frac{1}{|H|}\sum_{h\in H}\frac{\psi(h)}{|C_H(h)|} 
        & \text{if $\theta = 1_P$} \\ 
      0 & \text{if $\theta \neq 1_P$}  \\ 
    \end{array}\right. $$
In both cases, we have $[\chi,u] \in \CO$.  

Finally, we  remark that if $\CO$ is a Dedekind domain in which no 
prime dividing the order of $G$ is invertible, then  $A$ has the  
scalar projective property if and only if $G$ is abelian.  
\end{Example} 

\begin{Example}  \label{knorrmoritaexample}
The Kn\"orr property is not preserved by Morita equivalences in 
general. The idea is that all absolutely indecomposable $A$-lattices 
of $p'$-rank are Kn\"{o}rr, but among those of rank divisible by $p$, 
only the absolutely irreducible lattices tend to have the property. 
Indeed, the proof of \cite[Corollary 1.6]{Kno89} does not require the 
$\CO$-algebra to be a group ring (nor even a symmetric algebra). Thus, 
any Morita equivalence that sends a lattice of $p'$-rank which is 
indecomposable but not irreducible to a lattice of rank divisible by 
$p$ is likely to give an example.  

Specifically, let $p=2$ and assume that $\CO$ is unramified and  
$k$ is algebraically closed. Let $A$ be the principal block algebra  
of $\CO A_5$,  where $A_5$ is the alternating group of degree 5. Then 
$|\Irr_K(A)|=4 $, $|\Irr_k(A)| = 3$, and the decomposition matrix 
of $A$  with respect to some ordering of $\Irr_K(A) $  is
         \begin{equation}
                \begin{array}{cccc}
                      & \varphi_{1} & \varphi_{2}  &\varphi_3 \\\hline
                        \chi_{1} & 1 & 0 &0\\
                        \chi_{2} & 1 & 0 &1 \\ 
                        \chi_{3} & 1 & 1& 0 \\
                        \chi_{4} &1&1&1
                \end{array}
         \end{equation}
where $\varphi_ 1$ corresponds to a $1$-dimensional $kA$-module, and 
$\varphi_2 $ and $\varphi_3 $ correspond to simple $kA$-modules of   
dimension $2$. 

For each $i$, $1 \leq i \leq 3 $, let $P_i$ denote  a projective 
indecomposable  $A$-module  such that $P_i/\mathrm{rad}(P_i)$ is 
isomorphic to  a simple $kA$-module corresponding to $\varphi_i$.
 Let $e$ be an idempotent in $A$  such that 
$Ae \cong  P_1 + 2 P_2 + P_3 $ as left $A$-modules. Then $A$ and $eAe$ 
are Morita equivalent via the functor sending an $A$-module $M$ to the 
$eAe$-module $eM$   and an $A$-module   $\alpha: M \to N$ 
homomorphism  to the  $eAe$-module homomorphism  
$ e\cdot \alpha:  eM \to eN $   defined   through restriction to  
$eM$. The simple $eAe$-modules  corresponding to  $\varphi_1 $ and 
$\varphi_3 $  have  dimension $1$  whereas the simple $eAe$-module 
corresponding to $\varphi_2 $  has dimension $2$.

From the decomposition matrix above, one  sees that the character 
afforded by $KP_1$ has two irreducible constituents, one of degree 
$5$ and the other of degree $3$. 
It follows from \cite[Lemma 1.9]{Kno89} that $P_1$ is not a Kn\"{o}rr  
$A$-lattice.  However, the rank of the $eAe$-lattice $eP_1$ is $7$. 
Thus $eP_1$ is a Kn\"{o}rr $eAe$-lattice. 

To obtain an example in which neither lattice is projective, it is 
enough to inflate the $P_i$ above to lattices for the group 
$A_5\times C_2$, where $C_2$ is a cyclic group of order $2$.

Notice also that although $eP_1$ is Kn\"{o}rr, it does not have the 
stable exponent property. This is the case for both the $A_5$ and 
$A_5\times C_2$ situations. Next, we produce a lattice with the stable 
exponent property which is not Kn\"{o}rr.

First, we have $\Q(\sqrt{5}) \subseteq K,$  so $KA$ is split   
semisimple. Let $M$ be the unique quotient lattice of $P_1$  such that 
$KM $ has character $ \chi_1+ \chi_2 + \chi_3 $ .  Since  $M$ has rank 
$7$,  $M$ is a Kn\"{o}rr  $A$-lattice. Because $A$ has the projective 
scalar property, $M$ and hence $eM$ also have the stable exponent 
property. We shall show that $eM$ is not Kn\"{o}rr.

 Let $L$  be the unique $\CO$-free quotient of $M$ affording the 
character $\chi_3$ and let $\alpha: M \to L$ be the projection map. 
Since $\alpha$ is surjective,  and $L$ and  $M$ are not projective,    
$\alpha \notin \Hom_A^{\pr} ( M, L)$.
Thus, by Corollary \ref{TateDuality}, there exists   
$\beta  \in \Hom_A  (L, M) $ such that 
$\tr_M ( \beta \circ \alpha ) \notin 4 \CO $  
(since   $ 4 \cdot 1_A $ is a projective scalar element of $A$).     

Let $\tau =\beta \alpha $ and denote  also by $\tau $ the $K$-linear 
extension of $\tau$ to  $KM$.
For each $i$, $1\leq i \leq 4 $, let $e_i $ be the  primitive central 
idempotent of $KA$ corresponding to $\chi_i$. Since $\tau (K M)$ is 
contained in $e_3(K M)$, we have that $(e_2 +e_4)(K M)$ is contained 
in the kernel of $\tau$. On the other hand, $1-e = (e_2 +e_4) (1-e)$.
Thus, $\tr_{K M} (\tau) = \tr_{e(KM) }(\tau)$.   It follows that

 $$ \tr_{eM} (e\cdot \tau) =\tr_{K M} (\tau) = \tr_{M} (\tau)   
\notin 4\CO\ .  $$

Since   $eM$ has rank $6$,  we have that   
$\nu_2(\tr_{eM} (e\cdot \tau)) \leq \nu_2( \rank_{\CO} (eM))$. Since   
$\tau$ is not invertible,  neither is $e\cdot \tau$,  hence $eM$ is 
not a Kn\"{o}rr $eAe$-lattice.
\end{Example} 

\begin{Example}\label{non-converseexample}
Let $\CO=\Z_3$, and consider the $\CO$-order $A=\CO S_3$, that is, 
the group ring of the symmetric group on three points. The 
decomposition matrix of $A$ is 
         \begin{equation}
                \begin{array}{ccc}
                      & \varphi_{(3)} & \varphi_{(2,1)}  \\\hline
                        \chi_{(3)} & 1 & 0 \\
                        \chi_{(2,1)} & 1 & 1 \\
                        \chi_{(1^3)} & 0 & 1
                \end{array}
         \end{equation}
Here we use the  standard indexing of ordinary  and modular 
irreducible  characters of symmetric groups via partitions.  Let 
$e_{(3)}$, $e_{(2,1)}$ and $e_{(1^3)}$  denote the primitive 
idempotents in $Z(KA)$. The inertial index of this block is two, and, 
according to \cite{Bessenrodt}, that means that this block has six 
isomorphism types of indecomposable lattices (one can also show this  
in an elementary way). It is also easy to enumerate those isomorphism 
types: there are two indecomposable projective lattices, which are 
non-irreducible. Then there is a unique lattice with character 
$\chi_{(3)}$ and a unique lattice with character $\chi_{(1^3)}$. 
Moreover there is a lattice with character $\chi_{(2,1)}$ whose top 
has Brauer character $\varphi_{(3)}$ and there is a lattice with 
character $\chi_{(2,1)}$ whose top has Brauer character 
$\varphi_{(2,1)}$ (those two lattices are the projective lattices over 
the order $Ae_{(2,1)}$). As there are but six lattices in total we 
know that there can be no further indecomposable lattices. In 
particular, all indecomposable lattices are either projective or 
absolutely irreducible. This implies that each algebra in the Morita 
equivalence class of $A$ has the property that Kn\"orr-lattices and 
absolutely indecomposable non-projective lattices with the stable 
endomorphism property coincide. Any algebra in the Morita equivalence 
class of $A$ which does not possess the projective scalar property 
will therefore provide a counterexample to the converse of Theorem 
\ref{KnoerrExponent}.

Choose $B$ in the Morita equivalence class of $A$ such that the Morita 
equivalence sends the simple module with character $\varphi_{(3)}$ to 
a one-dimensional module and the simple module with character 
$\varphi_{(2,1)}$ to a two-dimensional module. Note that 
         \begin{equation}
         \frac{1}{3} \cdot \left(\chi_{(3)}(-) + 
         2\cdot \chi_{(2,1)}(-) + \chi_{(1^3)}(-)\right)
         \end{equation}
is a symmetrising form for $A$, and therefore also for $B$ (with the 
characters replaced by the corresponding characters of $B$). It 
follows that 
         \begin{equation}
          z_B = 3\cdot \left(e_{(3)} + \frac{3}{2} 
         \cdot e_{(2,1)} + 2\cdot  e_{(1^3)}\right)
         \end{equation}
and this element is determined uniquely up to multiplications by units 
in $Z(A)=Z(B)$. But multiplication by units cannot turn the above 
element into a scalar, since it will leave the $3$-valuation of the 
coefficients of the idempotents $e_{(3)}$, $e_{(2,1)}$ and 
$e_{(1^3)}$ invariant. Hence $B$ does not possess the projective 
scalar property.
\end{Example}

%%%%%%%%%%%%%%%%%%%%%%%%%%%%%%%%%%%%%%%%%%%%%%%%%%%%%%%%%%%%%%%%%%%%%%%
%%%%%%%%%%%%%%%%%%%%%%%%%%%%%%%%%%%%%%%%%%%%%%%%%%%%%%%%%%%%%%%%%%
%{\it Acknowledgement.} The present work is partially funded by
%the EPSRC grant EP/M02525X/1. 

%%%%%%%%%%%%%%%%%%%%%%%%%%%%%%%%%%%%%%%%%%%%%%%%%%%%%%%%%%%%%%%%%%%%%%

\end{document}